\documentclass[12pt]{article}
\usepackage{a4,latexsym, amssymb}

\usepackage[mathscr]{eucal}

\def \epf {\hfill $\Box$\\ \par}

\newtheorem{theorem}{Theorem}
\newtheorem{lemma}{Lemma}
\newtheorem{proposition}{Proposition}

\newtheorem{Remark}{Remark}

\font\msbm=msbm10 at 12pt

\begin{document}

\title{New Results on Two Hypercube Coloring Problems
\thanks{This research is supported in part by the National Natural Science Foundation of China under the Grant 60872025 and 10990011. The material in this paper was presented at
The Fifth Shanghai Conference on Combinatorics, May 14--18, 2005
}}

\author{ Fang-Wei Fu\thanks{Fang-Wei Fu is with the
Chern Institute of Mathematics and the Key Laboratory
of Pure Mathematics and Combinatorics, Nankai University, Tianjin 300071, P. R. China.
E-mail: fwfu@nankai.edu.cn},
San Ling\thanks{San Ling is with
the Division of Mathematical Sciences, School of Physical and
Mathematical Sciences, Nanyang Technological University,
SPMS-MAS-03-01, 21 Nanyang Link, Singapore 637371, Republic of Singapore.
E-mail: lingsan@ntu.edu.sg}, and Chaoping Xing\thanks{Chaoping Xing
is with the Division of Mathematical Sciences, School of Physical and
Mathematical Sciences, Nanyang Technological University,
SPMS-MAS-03-01, 21 Nanyang Link, Singapore 637371, Republic of Singapore.
E-mail: xingcp@ntu.edu.sg} }

\date{}
\maketitle

\begin{abstract}
In this paper, we study the following two hypercube coloring
problems: Given $n$ and $d$, find the minimum number of colors,
denoted as ${\chi}'_{d}(n)$ (resp. ${\chi}_{d}(n)$), needed to color
the vertices of the $n$-cube such that any two vertices with Hamming
distance at most $d$ (resp. exactly $d$) have different colors.
These problems originally arose in the study of the scalability of
optical networks. Using methods in coding theory, we show that
${\chi}'_{4}(2^{r+1}-1) = 2^{2r+1}$, ${\chi}'_{5}(2^{r+1}) =
4^{r+1}$ for any odd number $r\geq3$, and give two upper bounds on
${\chi}_{d}(n)$. The first upper bound improves on that of Kim, Du
and Pardalos. The second upper bound improves on the first one for
small $n$. Furthermore, we derive an inequality on ${\chi}_{d}(n)$
and ${\chi}'_{d}(n)$.
\end{abstract}


\pagebreak

\section{Introduction}

Let $V_n$ be the $n$-dimensional vector space over the binary
field ${\mathbb{F}}_{2}=\{0,1\}$, i.e.,
$$V_n = \{(a_1, a_2, \cdots, a_n) \mid a_i =0,1\}. $$
For ${\bf a}= (a_1, a_2, \cdots, a_n)\in V_n$ and ${\bf b}= (b_1,
b_2, \cdots, b_n)\in V_n$, the Hamming distance $d_{H}({\bf a},
{\bf b})$ between ${\bf a}$ and ${\bf b}$ is the number of
coordinates in which they differ, i.e.,
$$ d_{H}({\bf a}, {\bf b})= \mid \{i \mid a_i \not= b_i\}\mid.$$
The Hamming weight $w_{H}({\bf a})$ of a vector ${\bf a}\in V_n$
is the number of nonzero coordinates in ${\bf a}$. Obviously,
$$ d_{H}({\bf a}, {\bf b})= w_{H}({\bf a} - {\bf b}).$$

An $n$-cube (or $n$-dimensional hypercube) is a graph with the
vertex set $V_n$ and the edge set
$$ E_n = \{({\bf a}, {\bf b}) \mid {\bf a}, {\bf b} \in V_n,
d_{H}({\bf a}, {\bf b})=1 \}.$$ A coloring of $V_{n}$ with $L$
colors is a map $\Gamma$ from the vertex set $V_{n}$ to ${\cal L} =
\{1, 2, \cdots, L\}$. We say that the vertex ${\bf u}\in V_n$ is
colored with the color ${\Gamma}({\bf u})\in {\cal L}$. A
$d$-distance (resp. An exactly $d$-distance) coloring of $V_n$ is to
color the vertices of $V_n$ such that any two vertices with Hamming
distance at most $d$ (resp. exactly $d$) have different colors.
Note that for a coloring of $V_n$ with $L$ colors
\[ {\Gamma}: V_n \longrightarrow {\cal L}= \{1, 2, \cdots, L\},\]
it is a $d$-distance coloring of $V_n$ if and only if for any two
distinct vertices ${\bf u}, {\bf v} \in V_n$,
\[ {\Gamma}({\bf u}) \not= {\Gamma}({\bf v})
\;\;\;\mbox{if}\;\; d_{H}({\bf u}, {\bf v})\leq d;   \]
it is an
exactly $d$-distance coloring of $V_n$ if and only if for any two
distinct vertices ${\bf u}, {\bf v} \in V_n$,
\[ {\Gamma}({\bf u}) \not= {\Gamma}({\bf v})
\;\;\;\mbox{if}\;\; d_{H}({\bf u}, {\bf v})= d.   \] Denote
${\chi}'_{d}(n)$ (resp. ${\chi}_{d}(n)$) as the minimum number of
colors of a $d$-distance (resp. an exactly $d$-distance) coloring
of $V_n$.

These two hypercube coloring problems originally arose in the
study of the scalability of optical networks (see \cite{pav}). In
\cite{eno}-\cite{fran}, \cite{jm1}-\cite{lmt},
\cite{ngo1}-\cite{wan}, \cite{zhs} and \cite{zi}, some bounds on
${\chi}'_{d}(n)$ and ${\chi}_{d}(n)$ are given and some exact
values are determined.

The $d$-distance coloring and exactly $d$-distance coloring of $V_n$
are equivalent to certain partitions of $V_n$, which are related to
codes in coding theory. A nonempty subset $C$ of $V_n$ is called a
binary code of length $n$. Any word in $C$ is called a
codeword of $C$. The minimum distance $d(C)$ is defined as the
minimum distance between two distinct codewords of $C$. A binary
code $C$ of length $n$ and minimum distance $d$ is called a binary
$(n, d)$ code. A binary
code $C$ of length $n$ and minimum distance at least $d$ is called a binary
$(n, \geq d)$ code. For fixed length $n$ and minimum distance $d$, let
$A(n,d)$ denote the maximum size of a binary $(n, d)$ code. A binary
code $C$ is called a binary $[n,k]$ linear code if $C$ is a
$k$-dimensional subspace of $V_n$. Note that the size $|C|=2^k$ for
a binary $[n,k]$ linear code $C$. It is also well known in coding
theory that the minimum distance $d(C)$ of a linear code $C$ is
equal to the minimum Hamming weight of nonzero codewords. A binary
$[n,k]$ linear code with minimum distance $d$ is called a binary
$[n,k,d]$ linear code. For fixed length $n$ and minimum distance
$d$, let $k(n,d)$ denote the maximum dimension of a binary $[n,k,d]$
linear code. To determine $A(n,d)$ and $k(n,d)$ are two of the
fundamental research problems in coding theory. A binary code $C$ of
length $n$ is called an $(n, \overline{\{d\}})$ forbidden distance
code if $d_{H}({\bf a}, {\bf b})\not= d$ for any two distinct
codewords ${\bf a}, {\bf b} \in C$. Given $n$ and $d$, let
$Q(n,\overline{d})$ denote the maximum size of a binary $(n,
\overline{\{d\}})$ forbidden distance code.

An $(n,L,d)$-partition (resp. $(n,L,\overline{\{d\}})$-partition)
of $V_n$ is a set of subsets $\{B_i\}_{i=1}^{L}$ of $V_n$
satisfying {\rm (i)} $B_i \cap B_j = \emptyset$ for $i \not= j$
and ${\bigcup}_{i=1}^{L}B_i = V_n$; {\rm (ii)} Each $B_i$ is a
binary $(n, \geq d)$ code (resp. $(n, \overline{\{d\}})$ forbidden
distance code). It is well known that a $d$-distance coloring
(resp. an exactly $d$-distance coloring) of $V_n$ with $L$ colors
is equivalent to an $(n,L,d+1)$-partition (resp.
$(n,L,\overline{\{d\}})$-partition) of $V_n$ (note that
$B_i={\Gamma}^{-1}(i)$). Hence, ${\chi}'_{d}(n)$ (resp.
${\chi}_{d}(n)$) is the minimum number $L$ of subsets in any
$(n,L,d+1)$-partition (resp. $(n,L,\overline{\{d\}})$-partition)
of $V_n$.

Ngo et al. \cite{ngo1} obtained the following lower bound:
\begin{eqnarray}
{\chi}'_{d}(n) \geq \frac{2^n}{A(n,d+1)}. \label{ab1}
\end{eqnarray}
Similarly, for an $(n,L,\overline{\{d\}})$-partition $\{B_i\}_{i=1}^{L}$ of $V_n$,
$$ 2^n =|V_n| = \sum_{i=1}^{L}|B_i| \leq L\cdot Q(n,\overline{d}).  $$
Hence, we have
\begin{eqnarray}
{\chi}_{d}(n) \geq \frac{2^n}{Q(n,\overline{d})}. \label{ab2}
\end{eqnarray}
Note that for a binary $[n,k,d+1]$ linear code $C$, the cosets of
$C$ form an $(n,2^{n-k},d+1)$-partition of $V_n$. Hence, if there
exists a binary $[n,k,d+1]$ linear code, then ${\chi}'_{d}(n) \leq
2^{n-k}$. In particular, ${\chi}'_{d}(n) \leq 2^{n-k(n,d+1)}$.
Furthermore, if $A(n,d+1)=2^{k(n,d+1)}$, i.e., $A(n,d+1)$ is
attained by a binary $[n,k,d+1]$ linear code, then by (\ref{ab1}),
${\chi}'_{d}(n) = 2^{n-k(n,d+1)}$. In this way, a number of exact
values of ${\chi}'_{d}(n)$ are determined.

Our main results in this paper are given by the following theorems.
\begin{theorem}
For any odd integer $r\geq 3$, we have
\begin{eqnarray}
 {\chi}'_{4}(2^{r+1}-1) =2^{2r+1},
\label{aa1} \\
{\chi}'_{5}(2^{r+1}) = 4^{r+1}. \label{aa2}
\end{eqnarray}
\label{theorem1}
\end{theorem}
\begin{theorem}
For any even number $d$, we have
\begin{eqnarray}
{\chi}_{d}(n) \leq 2^{\left\lceil \log_{2}\left[1+ {n-1 \choose d-1}
\right]\right\rceil}  \label{aa3}
\end{eqnarray}
where $\lceil x \rceil$ is the smallest integer not less than $x$.
\label{theorem2}
\end{theorem}
\begin{theorem}
If $d$ is even and $n\geq 2d$, then
\begin{eqnarray}
{\chi}_{d}(n) \leq 2^{n-d+1- k(n-d+1, d+1)}.\label{bou}
\end{eqnarray}
\label{theorem3}
\end{theorem}

Theorem \ref{theorem1} could not be obtained by using the above
argument of binary linear codes. We prove Theorem \ref{theorem1}
in Section 2 by using the argument of ${\mathbb{Z}}_{4}$-linear
codes. In Section 3, we prove Theorems \ref{theorem2} and
\ref{theorem3} by using the argument of linear forbidden distance
codes. In order to prove Theorems \ref{theorem2} and
\ref{theorem3}, some properties of linear forbidden distance codes
are given. Furthermore, we derive the following inequality on
${\chi}_{d}(n)$ and ${\chi}'_{d}(n)$.
\begin{proposition}
\begin{eqnarray}
{\chi}_{d}(n_1+n_2) \leq {\chi}'_{d}(n_1) {\chi}_{d}(n_2), \;\;\;
\mbox{for even}\; d. \label{yy1}
\end{eqnarray}
\label{prop}
\end{proposition}

Kim, Du and Pardalos \cite{kim} showed that ${\chi}_{d}(n) = 2$
for odd $d$ and
\begin{eqnarray}
{\chi}_{d}(n) \leq
(d+1)\left(\frac{d+2}{2}\right)^{\frac{d(d+2)}{8}\lceil
{\log}_{2}n\rceil}, \;\; \mbox{for even} \; d. \label{aa4}
\end{eqnarray}
It is not hard to check that Bound (\ref{aa3}) improves on Bound
(\ref{aa4}). Note that Bound (\ref{bou}) depends on $k(n,d)$ and
is not explicit. Tables of lower and upper bounds on $k(n,d)$ for
small $n$ are available in \cite{brou} and \cite{brouw}. Using
Theorem \ref{theorem3} and presently best known lower bounds on
$k(n,d)$ listed in \cite{brou} and \cite{brouw}, we can obtain
upper bounds on ${\chi}_{d}(n)$ for small $n$, which improve on
Bound (\ref{aa3}). For example, if $n=13, d=4$, Bound (\ref{aa3})
gives ${\chi}_{4}(13) \leq 2^{8}$, while Bound (\ref{bou}) gives
${\chi}_{4}(13) \leq 2^{7}$; if $n=14, d=4$, Bound (\ref{aa3})
gives ${\chi}_{4}(14) \leq 2^{9}$, while Bound (\ref{bou}) gives
${\chi}_{4}(13) \leq 2^{7}$; if $n=28, d=6$, Bound (\ref{aa3})
gives ${\chi}_{6}(28) \leq 2^{17}$, while Bound (\ref{bou}) gives
${\chi}_{6}(28) \leq 2^{11}$.

\section{Proof of Theorem \ref{theorem1}}

In this section, we prove Theorem \ref{theorem1} by using the
argument of ${\mathbb{Z}}_{4}$-linear codes. Actually, Ngo et al.
\cite{ngo2} observed that ${\mathbb{Z}}_{4}$-linear codes may be
used to study this hypercube coloring problem. Unfortunately, they
\cite{ngo2} said that the construction of
${\mathbb{Z}}_{4}$-linear codes is quite involved, and gave no
results in this direction.

In order to establish our results, we first review some basic
definitions, notations and results for ${\mathbb{Z}}_{4}$-linear codes.
A detailed treatment of ${\mathbb{Z}}_{4}$-linear codes can be found in \cite{hamm} and \cite{wanz}.

Let ${\mathbb{Z}}_{4}$ be the ring of integers modulo $4$, and
let ${\mathbb{Z}}_{4}^{n}$ be the set of $n$-tuples over
${\mathbb{Z}}_{4}$. A nonempty subset ${\cal C}$ of ${\mathbb{Z}}_{4}^{n}$
is called a ${\mathbb{Z}}_{4}$-code. An additive subgroup ${\cal C}$ of ${\mathbb{Z}}_{4}^{n}$ is called a ${\mathbb{Z}}_{4}$-linear code. Any
${\mathbb{Z}}_{4}$-linear code is permutation equivalent to a
${\mathbb{Z}}_{4}$-linear code with a generator matrix of the
form
$$ G=\left(\begin{array}{ccc}
I_{k_1} & A        & B \\
O       & 2I_{k_2} & 2C
\end{array}\right) $$
where $A$ and $C$ are matrices with entries from $\{0,1\}$ and $B$
is a ${\mathbb{Z}}_{4}$-matrix. Then ${\cal C}$ is called a
${\mathbb{Z}}_{4}$-linear code of type $4^{k_1}2^{k_2}$, and the
size of ${\cal C}$ is $|{\cal C}|=4^{k_1}2^{k_2}$.

The Lee weights of $0, 1, 2, 3\in {\mathbb{Z}}_{4}$ are
defined by
$$ w_{L}(0)=0, \; w_{L}(1)=1, \; w_{L}(2)=2, \; w_{L}(3)=1. $$
The Lee weight $w_{L}({\bf x})$ of ${\bf x}=(x_1, \cdots, x_n) \in
{\mathbb{Z}}_{4}^{n}$ is defined to be the rational sum of the
Lee weights of its components, i.e.,
$$w_{L}({\bf x}) = \sum_{i=1}^{n} w_{L}(x_i).  $$
The Lee distance $d_{L}({\bf x}, {\bf y})$ between two tuples
${\bf x}, {\bf y} \in {\mathbb{Z}}_{4}^{n}$ is defined by
$$ d_{L}({\bf x}, {\bf y}) = w_{L}({\bf x}-{\bf y}).$$
For a ${\mathbb{Z}}_{4}$-code ${\cal C}$, the minimum Lee
distance $d_{L}({\cal C})$ of ${\cal C}$ is defined as the minimum
Lee distance between two distinct codewords of ${\cal C}$. It is
well known that for a ${\mathbb{Z}}_{4}$-linear code ${\cal
C}$, the minimum Lee distance $d_{L}({\cal C})$ is equal to the
minimum Lee weight $w_{L}({\cal C})$ of the nonzero codewords of
${\cal C}$.

The Gray map $\phi$ is a map from ${\mathbb{Z}}_{4}$ to
${\mathbb{Z}}_{2}^{2}$ defined as follows:
$$ {\phi}(0)=00, \; {\phi}(1)=01, \; {\phi}(2)=11, \; {\phi}(3)=10. $$
Clearly, $\phi$ is a bijection from ${\mathbb{Z}}_{4}$ to
${\mathbb{Z}}_{2}^{2}$. The Gray map $\phi$ can be extended in
an obvious way to a map from ${\mathbb{Z}}_{4}^{n}$ to
${\mathbb{Z}}_{2}^{2n}$ as follows:
$$ {\phi}({\bf x})= ({\phi}(x_1), \cdots, {\phi}(x_n)), \;\;
\mbox{for} \; {\bf x}=(x_1, \cdots, x_n) \in {\mathbb{Z}}_{4}^{n}. $$ Note that we still use $\phi$ to denote this
extended map. It is easy to see that $\phi$ is still a bijection
from ${\mathbb{Z}}_{4}^{n}$ to ${\mathbb{Z}}_{2}^{2n}$.
The following result about the Gray map can be found in
\cite{hamm} and \cite{wanz}.
\begin{lemma}
The Gray map $\phi$ is a distance-preserving map from $({\mathbb{Z}}_{4}^{n},\; d_{L})$
to $({\mathbb{Z}}_{2}^{2n},\; d_{H})$, i.e., $d_{L}({\bf x},
{\bf y}) = d_{H}({\phi}({\bf x}), {\phi}({\bf y}))$ for all ${\bf
x}, {\bf y} \in {\mathbb{Z}}_{4}^{n}$, and $\phi$ is also a
weight-preserving map from $({\mathbb{Z}}_{4}^{n},\; w_{L})$ to $({\mathbb{Z}}_{2}^{2n},\; w_{H})$,
i.e., $ w_{L}({\bf x}) = w_{H}({\phi}({\bf x}))$ for all ${\bf x}
\in {\mathbb{Z}}_{4}^{n}$. \label{lemma1}
\end{lemma}
Let ${\cal C}$ be a ${\mathbb{Z}}_{4}$-code of length $n$. The
binary image of ${\cal C}$ is defined as
$$ C = {\phi}({\cal C}) = \{{\phi}({\bf c})\;|\; {\bf c} \in {\cal C}\}.  $$
Clearly, $C$ is a binary code of length $2n$ and size $|C|=|{\cal
C}|$. It follows from Lemma \ref{lemma1} that the minimum Hamming
distance of $C$ is equal to the minimum Lee distance of ${\cal C}$,
i.e., $d(C) = d_{L}({\cal C})$. Note that the binary image of a
${\mathbb{Z}}_{4}$-linear code may not be a binary linear code.
Some well-known nonlinear binary codes such as the
Nordstrom-Robinson, Kerdock, Preparata, Goethals, and
Delsarte-Goethals codes contain more codewords than any known binary
linear codes with the same length and minimum distance. For example,
the Preparata code contains twice as many codewords as the extended
$2$-error-correcting BCH code of the same length. Hammons et al.
\cite{hamm} showed that these codes can be very simply constructed
as the binary image of certain ${\mathbb{Z}}_{4}$-linear codes
under the Gray map.

In the following we use ${\mathbb{Z}}_{4}$-linear codes to
construct $(n,L,d)$-partitions of $V_n$.
\begin{theorem}
If there exists a ${\mathbb{Z}}_{4}$-linear code of length $n$,
type $4^{k_1}2^{k_2}$ and minimum Lee distance $d_{L}\geq 3$, then
there exist a $(2n, 2^{2n-2k_{1}-k_{2}}, d_{L})$-partition of
$V_{2n}$ and a $(2n-1, 2^{2n-1-2k_{1}-k_{2}}, d_{L}-1)$-partition of
$V_{2n-1}$.
\label{theorem4}
\end{theorem}
\begin{proof}
Let ${\cal C}$ be a ${\mathbb{Z}}_{4}$-linear
code of length $n$, type $4^{k_1}2^{k_2}$ and minimum Lee distance
$d_{L}\geq 3$. Recall that the size $|{\cal C}|=4^{k_1}2^{k_2}$
and $d_{L}$ is equal to the minimum Lee weight of ${\cal C}$. Note
that ${\cal C}$ is an additive subgroup of ${\mathbb{Z}}_{4}^{n}$. The number of cosets of ${\cal C}$ is given by
$$ N = \frac{4^n}{4^{k_1}2^{k_2}} = 2^{2n-2k_{1}-k_{2}}. $$
Denote these cosets as
$$ {\cal B}_{i} = {\bf b}_{i} + {\cal C}, \;\; i=1, 2, \cdots, N  $$
where ${\bf b}_{1} = {\bf 0}$ and ${\bf b}_{i} \in {\mathbb{Z}}_{4}^{n}$.
Note that every ${\cal B}_{i}$ is a ${\mathbb{Z}}_{4}$-code of length $n$, size $4^{k_1}2^{k_2}$ and minimum Lee
distance $d_{L}$. Furthermore, $\{{\cal B}_{i}\}_{i=1}^{N}$ forms
a partition of ${\mathbb{Z}}_{4}^{n}$. Denote the binary image
of ${\cal B}_{i}$ under the Gray map $\phi$ as
$$ B_{i} = {\phi}({\cal B}_{i}), \;\; i=1, 2, \cdots, N. $$
Recall that $\phi$ is a bijection from ${\mathbb{Z}}_{4}^{n}$ to
${\mathbb{Z}}_{2}^{2n}$. By Lemma \ref{lemma1}, we know that
$B_i$ is a binary code of length $2n$, size $4^{k_1}2^{k_2}$ and
minimum Hamming distance $d_{L}$. Furthermore, $\{B_{i}\}_{i=1}^{N}$
forms a partition of ${\mathbb{Z}}_{2}^{2n}$. Hence,
$\{B_{i}\}_{i=1}^{N}$ forms a $(2n, N, d_{L})$-partition of
$V_{2n}$.

Next, we further construct a $(2n-1, N/2, d_{L}-1)$-partition of
$V_{2n-1}$. Consider the cosets of ${\cal C}$ with respect to the
additive group ${\mathbb{Z}}_{4}^{n}$, i.e., ${\cal B}_{1},
{\cal B}_{2}, \cdots, {\cal B}_{N}$. We say that two cosets ${\cal
B}_{i}$ and ${\cal B}_{j}$ are equivalent if there exists $s\in
{\mathbb{Z}}_{4}$ such that
\[ {\cal B}_{j} = (0, \cdots, 0, s) + {\cal B}_{i}. \]
Hence, the equivalence class containing the coset ${\cal B}_{i}$
is given as
\[ \widetilde{{\cal B}}_{i} = \{ {\cal B}_{i},\; (0, \cdots, 0, 1) + {\cal B}_{i},\;
(0, \cdots, 0, 2) + {\cal B}_{i},\; (0, \cdots, 0, 3) + {\cal B}_{i}
\}. \] The four cosets in the equivalence class $\widetilde{{\cal
B}}_{i}$ are distinct since $d_{L}\geq 3$. Let $N_{1}=N/4$. Note
that there are $N_{1}$ equivalence classes in total, and these
$N_{1}$ equivalence classes form a partition of the coset set
$\{{\cal B}_{1}, {\cal B}_{2}, \cdots, {\cal B}_{N}\}$. Denote these
$N_1$ equivalence classes as
\[\widetilde{{\cal B}}_{j_1}, \widetilde{{\cal B}}_{j_{2}}, \cdots, \widetilde{{\cal B}}_{j_{N_1}}. \]
Next we consider the binary images of four cosets in an equivalence
class. Hammons et al. \cite{hamm} observed that for ${\bf x}=(x_1,
\cdots, x_n) \in {\mathbb{Z}}_{4}^{n}$ and ${\bf y}=(y_1,
\cdots, y_n) \in {\mathbb{Z}}_{4}^{n}$
\[ {\phi}({\bf x}+ {\bf y})= {\phi}({\bf x}) + {\phi}({\bf y})
+ {\phi}(2{\alpha}({\bf x})*{\alpha}({\bf y}) )   \]
where
$\alpha$ is a map from ${\mathbb{Z}}_{4}$ to ${\mathbb{Z}}_{2}$ defined by
\[ {\alpha}(0)=0,\; {\alpha}(1)=1,\; {\alpha}(2)=0,\; {\alpha}(3)=1, \]
${\alpha}({\bf x})$ is defined as ${\alpha}({\bf x})=
({\alpha}(x_1), \cdots, {\alpha}(x_n))$ and ${\alpha}({\bf
x})*{\alpha}({\bf y})$ is defined as
\[ {\alpha}({\bf x})*{\alpha}({\bf y}) =
({\alpha}(x_1){\alpha}(y_1), \cdots, {\alpha}(x_n){\alpha}(y_n)).
\] If ${\bf y}=(0, \cdots, 0, 1) \in {\mathbb{Z}}_{4}^{n}$,
then
\[ {\phi}({\bf y}) = (0,0,\cdots,0,0,0,1) \in V_{2n} \]
and
\[ {\alpha}({\bf x})*{\alpha}({\bf y}) =
(0,\cdots,0,{\alpha}(x_n)) \in V_{n}. \] Hence, if
${\alpha}(x_n)=0$,
\[ {\phi}({\bf x}+ {\bf y})= {\phi}({\bf x}) + (0,0,\cdots,0,0,0,1),\]
if ${\alpha}(x_n)=1$,
\begin{eqnarray*}
{\phi}({\bf x}+ {\bf y})
& = & {\phi}({\bf x}) + (0,0,\cdots,0,0,0,1) + (0,0,\cdots,0,0,1,1) \\
& = & {\phi}({\bf x}) + (0,0,\cdots,0,0,1,0).
\end{eqnarray*}
If ${\bf y}=(0, \cdots, 0, 2) \in {\mathbb{Z}}_{4}^{n}$, then
\[ {\phi}({\bf y}) = (0,0,\cdots,0,0,1,1) \in V_{2n} \]
and
\[ {\alpha}({\bf x})*{\alpha}({\bf y}) =
(0,\cdots,0,0) \in V_{n}. \] Hence,
\[ {\phi}({\bf x}+ {\bf y})= {\phi}({\bf x}) + (0,0,\cdots,0,0,1,1).\]
If ${\bf y}=(0, \cdots, 0, 3) \in {\mathbb{Z}}_{4}^{n}$, then
\[ {\phi}({\bf y}) = (0,0,\cdots,0,0,1,0) \in V_{2n} \]
and
\[ {\alpha}({\bf x})*{\alpha}({\bf y}) =
(0,\cdots,0,{\alpha}(x_n)) \in V_{n}. \] Hence, if
${\alpha}(x_n)=0$,
\[ {\phi}({\bf x}+ {\bf y})= {\phi}({\bf x}) + (0,0,\cdots,0,0,1,0),\]
if ${\alpha}(x_n)=1$,
\begin{eqnarray*}
{\phi}({\bf x}+ {\bf y})
& = & {\phi}({\bf x}) + (0,0,\cdots,0,0,1,0) + (0,0,\cdots,0,0,1,1) \\
& = & {\phi}({\bf x}) + (0,0,\cdots,0,0,0,1).
\end{eqnarray*}
Therefore, for the four cosets in an equivalence class
$\widetilde{{\cal B}}_{i}$, the binary images are given by
\[ B_{i} = {\phi}({\cal B}_{i}) = \{{\phi}({\bf x})\;|\; {\bf x} \in {\cal B}_{i}\}, \]
\begin{eqnarray*}
B_{i,1} & = & {\phi}({\cal B}_{i}+(0, \cdots, 0, 1)) \\
& = & \left\{{\phi}({\bf x}) + (0,0,\cdots,0,0,0,1)\;|\; {\bf x}
\in {\cal B}_{i}\;
\mbox{and}\; {\alpha}(x_n)=0 \right\} \bigcup \\
& & \left\{{\phi}({\bf x}) + (0,0,\cdots,0,0,1,0)\;|\; {\bf x} \in
{\cal B}_{i}\; \mbox{and}\; {\alpha}(x_n)=1 \right\},
\end{eqnarray*}
\begin{eqnarray*}
B_{i,2} & = & {\phi}({\cal B}_{i}+(0, \cdots, 0, 2))
 = \left\{{\phi}({\bf x}) + (0,0,\cdots,0,0,1,1)\;|\; {\bf x} \in {\cal B}_{i}\right\},
\end{eqnarray*}
\begin{eqnarray*}
B_{i,3} & = & {\phi}({\cal B}_{i}+(0, \cdots, 0, 3)) \\
& = & \left\{{\phi}({\bf x}) + (0,0,\cdots,0,0,1,0)\;|\; {\bf x}
\in {\cal B}_{i}\;
\mbox{and}\; {\alpha}(x_n)=0 \right\} \bigcup \\
& & \left\{{\phi}({\bf x}) + (0,0,\cdots,0,0,0,1)\;|\; {\bf x} \in
{\cal B}_{i}\; \mbox{and}\; {\alpha}(x_n)=1 \right\}.
\end{eqnarray*}
For a nonempty subset $B$ of $V_{2n}$, denote by $B^{*}$ the set
obtained by deleting the last coordinate of each word in $B$. Since
$B_{i}$, $B_{i,1}$, $B_{i,2}$, $B_{i,3}$ are binary codes of length
$2n$, size $4^{k_1}2^{k_2}$ and minimum Hamming distance $d_{L}$,
then $B_{i}^{*}$, $B_{i,1}^{*}$, $B_{i,2}^{*}$, $B_{i,3}^{*}$ are
binary codes of length $2n-1$, size $4^{k_1}2^{k_2}$ and minimum
Hamming distance at least $d_{L}-1$. Furthermore, by the expressions
of $B_{i}$, $B_{i,1}$, $B_{i,2}$, $B_{i,3}$, we have
\[ B_{i}^{*} \cup B_{i,2}^{*} = B_{i,1}^{*} \cup B_{i,3}^{*}. \]
Applying the above procedure to every equivalence class of
\[\widetilde{{\cal B}}_{j_1}, \widetilde{{\cal B}}_{j_{2}}, \cdots, \widetilde{{\cal B}}_{j_{N_1}}, \]
we know that
\[ B_{j_1}^{*} \bigcup  B_{j_{1},2}^{*} \bigcup B_{j_2}^{*} \bigcup  B_{j_{2},2}^{*} \bigcup
\cdots \bigcup B_{j_{N_{1}}}^{*}
\bigcup B_{j_{N_{1}},2}^{*} = V_{2n-1}. \] Since
\[ |B_{j_1}^{*}| + |B_{j_{1},2}^{*}| + |B_{j_2}^{*}| + |B_{j_{2},2}^{*}| + \cdots + |B_{j_{N_{1}}}^{*}|
+ |B_{j_{N_{1}},2}^{*}| = 2N_{1}4^{k_1}2^{k_2} = 2^{2n-1}
=|V_{2n-1}|, \] then
\[ B_{j_1}^{*}, B_{j_{1},2}^{*}, B_{j_2}^{*}, B_{j_{2},2}^{*}, \cdots, B_{j_{N_{1}}}^{*}, B_{j_{N_{1}},2}^{*}\]
form a $(2n-1, N/2, d_{L}-1)$-partition of $V_{2n-1}$.
\end{proof}

Now we prove Theorem \ref{theorem1} by using Theorem \ref{theorem4}.

{\bf Proof of Theorem \ref{theorem1}}: For odd $r\geq 3$, Hammons et
al. \cite{hamm} constructed a ${\mathbb{Z}}_{4}$-linear code
${\cal P}(r)$ of length $2^{r}$, type $4^{2^{r}-r-1}$ and minimum
Lee distance $6$. The binary image $P(r)= {\phi}\left({\cal
P}(r)\right)$ is a nonlinear binary code of length $2^{r+1}$, size
$2^{2^{r+1}-2r-2}$ and minimum Hamming distance $6$. This binary
code is also called the Preparata code in \cite{hamm}. By Theorem \ref{theorem4}, there exist a $(2^{r+1}, 2^{2r+2}, 6)$-partition of
$V_{2^{r+1}}$ and a $(2^{r+1}-1, 2^{2r+1}, 5)$-partition of
$V_{2^{r+1}-1}$ for odd $r\geq 3$. Since a $d$-distance coloring of
$V_n$ with $L$ colors is equivalent to an $(n,L,d+1)$-partition of
$V_n$, we have for odd $r\geq 3$,
\begin{eqnarray}
 {\chi}'_{4}(2^{r+1}-1) \leq 2^{2r+1},
\;\;\;\; {\chi}'_{5}(2^{r+1}) \leq 4^{r+1}. \label{bb1}
\end{eqnarray}
Furthermore, the Preparata code is an optimal binary code, i.e.,
$A(2^{r+1},6)=2^{2^{r+1}-2r-2}$ for odd $r\geq 3$. It follows from
(\ref{ab1}) that
\begin{eqnarray}
 {\chi}'_{5}(2^{r+1}) \geq 4^{r+1}. \label{bb2}
\end{eqnarray}
Since $A(n,d)=A(n+1, d+1)$ for odd $d$ (see \cite{ms}), we have
$A(2^{r+1}-1,5)=A(2^{r+1},6)=2^{2^{r+1}-2r-2}$ for odd $r\geq 3$.
Using (\ref{ab1}) again, we obtain
\begin{eqnarray}
 {\chi}'_{4}(2^{r+1}-1) \geq 2^{2r+1}.
\label{bb3}
\end{eqnarray}
Theorem \ref{theorem1} follows from combining (\ref{bb1}) with
(\ref{bb2}) and (\ref{bb3}). \epf

\begin{Remark} By generalizing the above argument, if
there exists a ${\mathbb{Z}}_{4}$-linear code of length $n$, type
$4^{k_1}2^{k_2}$ and minimum Lee distance $d_{L} \geq 2t+1$, we
can construct a $(2n-t, N/2^{t}, d_{L}-t)$-partition of
$V_{2n-t}$. The proof of this fact is more tedious and omitted
here. In this way, we can construct a $(d_{L}-t-1)$-distance coloring of $V_{2n-t}$ with $N/2^{t}$ colors.
In some cases, this method yields good distance colorings of hypercubes with less colors.
\end{Remark}

\section{Proofs of Theorems \ref{theorem2} and \ref{theorem3} and Proposition \ref{prop}}

In this section, we prove Theorems \ref{theorem2} and
\ref{theorem3} and Proposition \ref{prop} by using the argument of
linear forbidden distance codes. Some properties of linear
forbidden distance codes are given.

A binary $[n,k]$ linear code $C$ is called a binary $[n,k,\overline{\{d\}}]$ linear
forbidden distance code if $C$ is also an $(n, \overline{\{d\}})$
forbidden distance code, that is, $w_{H}({\bf c})\not=d$ for any nonzero codeword of $C$. Given $n$ and $d$, denote by $k(n,
\overline{d})$ the maximum dimension of a binary
$[n,k,\overline{\{d\}}]$ linear forbidden distance code. The problem
of determining $k(n, \overline{d})$ is studied in \cite{bas} and
\cite{eno}. Note that for a binary $[n,k,\overline{\{d\}}]$ linear
forbidden distance code $C$, the cosets of $C$ form an
$(n,2^{n-k},\overline{\{d\}})$-partition of $V_n$. Recall that an
exactly $d$-distance coloring of $V_n$ with $L$ colors is equivalent
to an $(n,L,\overline{\{d\}})$-partition of $V_n$. Hence, we have
the following results.
\begin{lemma}
If there exists a binary $[n,k,\overline{\{d\}}]$ forbidden distance
linear code, then
\begin{eqnarray}
{\chi}_{d}(n) \leq 2^{n-k}. \label{ca1}
\end{eqnarray}
In particular,
\begin{eqnarray}
{\chi}_{d}(n) \leq 2^{n- k(n, \overline{d})}. \label{ca2}
\end{eqnarray}
\label{lemma3}
\end{lemma}

Below we derive some bounds on $k(n, \overline{d})$. Using these bounds and Lemma \ref{lemma3}, we prove Theorems \ref{theorem2} and \ref{theorem3}. For a binary matrix $H$ of size $m\times n$, we can define a linear
code
\begin{eqnarray}
C= \{{\bf c}\in V_{n}: H {\bf c}^{T} = {\bf 0}\}. \label{ca3}
\end{eqnarray}
The dimension of $C$, $dim(C)=n- r(H) \geq n-m$, where $r(H)$ is the
rank of $H$. It is not hard to see that $C$ is a binary
$(n,\overline{\{d\}})$ forbidden distance code if and only if the
sum of any $d$ columns of $H$ is not the zero vector, i.e., any
column of $H$ is not equal to the sum of any other $d-1$ columns of
$H$. Similar to the Varshamov-Gilbert bound (see \cite{ms}) for
linear codes in coding theory, we have the following lower bound on
$k(n, \overline{d})$.
\begin{lemma}
\begin{eqnarray}
k(n, \overline{d}) \geq n - \left\lceil \log_{2}\left[1+ {n-1
\choose d-1} \right]\right\rceil \label{ca4}
\end{eqnarray}
where $\lceil y \rceil$ is the smallest integer greater than or
equal to $y$. \label{lemma4}
\end{lemma}
\begin{proof}
We shall show that if
\begin{eqnarray}
1+ {n-1 \choose d-1} \leq 2^{m}, \label{ca5}
\end{eqnarray}
then there exists a binary matrix $H$ of size $m\times n$ such that
any column of $H$ is not equal to the sum of any other $d-1$ columns
of $H$. We choose the first column ${\bf h}_{1}$, $\cdots$, the
$n$-th column ${\bf h}_{n}$ of $H$ from $V_m$ step by step. In
general, we can choose ${\bf h}_{1}$ to be any vector of $V_m$, and
choose ${\bf h}_{j}$ to be any vector of $V_m$ that is not the sum
of any $d-1$ distinct vectors of ${\bf h}_{1}$, $\cdots$, ${\bf
h}_{j-1}$. Now, there are at most ${j-1 \choose d-1}$ such sum
vectors. It follows from (\ref{ca5}) that for $j=2, 3, \cdots, n$,
\[{j-1 \choose d-1} \leq {n-1 \choose d-1} \leq 2^{m}-1 < |V_m|.  \]
Hence, in this way we can choose the first column ${\bf h}_{1}$,
$\cdots$, the $n$-th column ${\bf h}_{n}$ of $H$. The binary linear
code $C$ defined by (\ref{ca3}) is a binary $[n,k,\overline{\{d\}}]$
linear forbidden distance code, where the dimension $k\geq n-m$. By
(\ref{ca5}), we can take
\[m= \left\lceil \log_{2}\left[1+ {n-1 \choose d-1} \right]\right\rceil.  \]
Hence,
\begin{eqnarray*}
k(n, \overline{d})
& \geq & dim(C) \\
& \geq & n-m \\
& =    & n - \left\lceil \log_{2}\left[1+ {n-1 \choose d-1}
\right]\right\rceil.
\end{eqnarray*}
This completes the proof.
\end{proof}

{\bf Proof of Theorem \ref{theorem2}}: Theorem \ref{theorem2}
follows from (\ref{ca2}) and (\ref{ca4}). \epf

We see from Lemma \ref{lemma3} that $k(n, \overline{d})$ can be
used to give an upper bound on ${\chi}_{d}(n)$. We next give some
properties of $k(n, \overline{d})$. Using these properties we
prove Theorem \ref{theorem3}.
\begin{proposition}
\begin{eqnarray}
k(n, \overline{2}) = n - \lceil \log_{2}n \rceil. \label{kk1}
\end{eqnarray}
\label{prop1}
\end{proposition}
\begin{proof} It follows from Lemma \ref{lemma4} that
\begin{eqnarray}
k(n, \overline{2}) \geq n - \lceil \log_{2}n \rceil. \label{ca6}
\end{eqnarray}
On the other hand, let $C$ be a binary $[n,k(n,
\overline{2}),\overline{\{2\}}]$ linear forbidden distance code.
For $i=1, 2, \cdots, n$, let ${\bf e}_{i}$ be the $i$-th unit
vector that only has a $1$ at the $i$-th position. Note that
$w_{H}({\bf c})\not=2$ for any nonzero codeword ${\bf c} \in C$.
Hence,
\[ \{{\bf e}_{i}+ C\} \cap \{{\bf e}_{j}+ C\} = \emptyset, \;\;\mbox{for} \; i\not= j \]
since ${\bf e}_{i}+{\bf e}_{j} \not\in C$. Obviously,
\[ \bigcup_{i=1}^{n} \{{\bf e}_{i}+ C\} \subseteq V_{n}.\]
Hence,
\[ 2^n = |V_n| \geq \sum_{i=1}^{n} |\{{\bf e}_{i}+ C\}| = n |C| = n 2^{k(n, \overline{2})}. \]
This implies that
\begin{eqnarray}
k(n, \overline{2}) \leq n - \lceil \log_{2}n \rceil. \label{ca7}
\end{eqnarray}
Combining (\ref{ca6}) and (\ref{ca7}), we obtain (\ref{kk1}). This
completes the proof.
\end{proof}
\begin{proposition}
\begin{eqnarray}
k(n+n', \overline{d})\geq k(n,d+1) + k(n', \overline{d}), \;\;\;
\mbox{for even}\; d. \label{kk2}
\end{eqnarray}
In particular,
\begin{eqnarray}
k(n+d-1, \overline{d})\geq k(n,d+1) + d-1,  \;\;\;\mbox{for
even}\; d. \label{kk3}
\end{eqnarray}
\label{prop2}
\end{proposition}
\begin{proof}  Let $C_1$ be a binary $[n,k(n, d+1),d+1]$ linear
code. Let $C_2$ be a binary $[n',k(n',
\overline{d}),\overline{\{d\}}]$ linear forbidden distance code.
Denote
\[ C= \{ ({\bf a}, {\bf b}): {\bf a}\in C_1, {\bf b}\in C_2\}.\]
It is easy to check that $C$ is a binary $[n+n',k(n,d+1)+k(n',
\overline{d}),\overline{\{d\}}]$ linear forbidden distance code.
Hence, by the definition of $k(n+n', \overline{d})$, we have
\begin{eqnarray*}
k(n+n', \overline{d})\geq k(n,d+1) + k(n', \overline{d}).
\end{eqnarray*}
Taking $n'=d-1$ in (\ref{kk2}) and noting that $k(d-1,
\overline{d})=d-1$, we obtain (\ref{kk3}).
\end{proof}

It follows from (\ref{kk3}) that $k(2d, \overline{d})\geq d$ for
even $d$. Enomoto {\sl et al}. \cite{eno} proved that $k(2d,
\overline{d})=d$ for even $d$. This shows that Bound (\ref{kk3})
is tight for some cases.

{\bf Proof of Theorem \ref{theorem3}}: It follows from (\ref{kk3})
that, if $d$ is even and $n\geq 2d$,
\begin{eqnarray}
k(n, \overline{d})\geq k(n-d+1,d+1) + d-1. \label{kk4}
\end{eqnarray}
Theorem \ref{theorem3} follows from (\ref{ca2}) and (\ref{kk4}).
\epf

Now we prove Proposition \ref{prop} by using the following lemma.
\begin{lemma}
If there exist an $(n_1,L_1,d+1)$-partition of $V_{n_1}$ and an
$(n_2,L_2,\overline{\{d\}})$-partition of $V_{n_2}$, then there
exists an $(n_1+n_2,L_1L_2,\overline{\{d\}})$-partition of
$V_{n_{1}+n_{2}}$. \label{lemma5}
\end{lemma}
\begin{proof} Suppose that the $L_1$ subsets $\{B_i\}_{i=1}^{L_1}$
of $V_{n_1}$ form an $(n_1,L_1,d+1)$-partition of $V_{n_{1}}$, and
the $L_2$ subsets $\{E_j\}_{j=1}^{L_2}$ of $V_{n_{2}}$ form an
$(n_2,L_2,\overline{\{d\}})$-partition of $V_{n_{2}}$. Then, by
the definitions of an $(n_1,L_1,d+1)$-partition of $V_{n_{1}}$ and
an $(n_2,L_2,\overline{\{d\}})$-partition of $V_{n_{2}}$, $B_i$ is
a binary $(n_1, \geq(d+1))$ code and $E_j$ is a binary $(n_2,
\overline{\{d\}})$ forbidden distance code. It is easy to see that
\[ B_i \times E_j = \{ ({\bf a}, {\bf b}): {\bf a}\in B_i, {\bf b}\in E_j\}\]
is an $(n_1+n_2, \overline{\{d\}})$ forbidden distance code.
Hence, the $L_1L_2$ subsets $B_i \times E_j$, $1\leq i\leq L_1$,
$1\leq j\leq L_2$, of $V_{n_{1}+n_{2}}$ form an
$(n_1+n_2,L_1L_2,\overline{\{d\}})$-partition of
$V_{n_{1}+n_{2}}$.
\end{proof}

{\bf Proof of Proposition \ref{prop}}: Note that a $d$-distance
coloring of $V_{n}$ with $L$ colors is equivalent to an
$(n,L,d+1)$-partition of $V_n$, and an exactly $d$-distance
coloring of $V_n$ with $L$ colors is equivalent to an
$(n,L,\overline{\{d\}})$-partition of $V_n$. Lemma \ref{lemma5}
tells us that for a $d$-distance coloring of $V_{n_1}$ with
$L_1={\chi}'_{d}(n_1)$ colors and an exactly $d$-distance coloring
of $V_{n_2}$ with $L_2={\chi}_{d}(n_2)$ colors, we can construct
an exactly $d$-distance coloring of $V_{n_{1}+n_{2}}$ with
$L_1L_2$ colors. Hence, we have
\begin{eqnarray*}
{\chi}_{d}(n_1+n_2) \leq {\chi}'_{d}(n_1) {\chi}_{d}(n_2), \;\;\;
\mbox{for even}\; d.
\end{eqnarray*}

\end{document}